\documentstyle[studies,named,twoside]{article}

\newfont{\bb}{msbm10}
\newfont{\sbm}{msbm10}
\newfont{\sam}{msam10}
\newcommand{\restriction}{\mbox{\sam \symbol{'026}}}
\newcommand{\nleq}{\mbox{\sbm \symbol{'004}}}
\newcommand{\Vdash}{\mbox{\sam \symbol{'015}}}
\newcommand{\nsubseteq}{\mbox{\sbm \symbol{'052}}}

\begin{document}
\bibliographystyle{named}

\title{Subalgebras of Cohen algebras need not be Cohen}
\ShortTitle{Subalgebras of Cohen algebras}                    
\author{Sabine Koppelberg and Saharon Shelah}
\ShortAuthor{S.\ Koppelberg and S.\ Shelah}
\maketitle



\addtocounter{section}{-1}
\section{Introduction}   

Let us denote by   $\mbox{\bb C}_\kappa $  the standard Cohen algebra of  
$\pi$-weight $\kappa$, i.e.\ the complete Boolean algebra adjoining  
$\kappa$  Cohen reals, where  $\kappa$  is an infinite cardinal or $0$. 
More 
generally, we call a Boolean algebra  $A$  a Cohen algebra if (for 
technical convenience in Theorems \ref{th:0.1} and \ref{th:0.2}  below)
it  satisfies the countable chain condition 
and forcing with  $A$  (more precisely with the partial ordering  
$A \setminus \{0\}$) is equivalent to Cohen forcing, i.e.\ if every
generic extension  of the universe 
of set theory arising from forcing with $A$  arises from forcing with 
some standard Cohen algebra.  Since forcing with an arbitrary Boolean
algebra is  equivalent to forcing with its completion and forcing with a
product of algebras is  equivalent to forcing with one of the factors,
an algebra  is Cohen iff its  completion is isomorphic to a product of
at most countably many standard  Cohen algebras; we will use this
description  as the definition of a Cohen algebra in the rest of the 
paper. 

Cohen algebras are among the most important objects to be studied in the 
realm of Boolean algebras or forcing. There is a general feeling that
more or less  ``every" algebraic property of the standard Cohen algebras
is well-known;  similarly, the effect of adding Cohen reals to a given
model of set theory is, generally, quite well understood. It is
therefore quite surprising that the answer to an  apparently innocent
question was open, up to now; cf.\ Problem 5.2 in 
\cite{Koppelberg:2}.

\begin{definition}[Problem]  If  $B$  is a regular subalgebra of some 
Cohen algebra  $A$, does it follow that  $B$  is Cohen?  
\end{definition}

The first reference to this problem we are aware of is in 
Kamburelis\rq s paper  \cite{Kamburelis:1}. By Example 5.5 in
\cite{Koppelberg:2}, the 
assumption 
that  $B$  be regular in  $A$  cannot be disposed with.

We make some simple observations which somewhat restrict the problem 
(cf.\ \cite{Koppelberg:2}). First, both  $A$  and  $B$  may be assumed to
be complete; 
moreover  $A$  may be assumed to be a standard Cohen algebra  $\mbox{\bb 
C}_\kappa $.
Finally,  $\kappa$  may be assumed to be at least  $\omega_2$, by 
Proposition 5.4 in  \cite{Koppelberg:2}. Thus the following theorem, 
the unique result of 
the paper, is the strongest result one can hope for.    

\begin{theorem} For every  $\kappa \geq  \omega_2$,  $\mbox{\bb 
C}_\kappa$  
has a complete regular subalgebra of  $\pi$-weight  $\kappa$ which is 
not Cohen. 
\end{theorem}

Let us mention that several beautiful internal descriptions of Cohen 
algebras were proved in recent years, based on 
Shapiro\rq s theorem (cf.\ \cite{Shapiro:1}, \cite{Shapiro:2}) that 
every subalgebra of a free  Boolean algebra is Cohen; see Section 1 for
all unexplained notions.  Two of these are given below; they will,
however, not be applied in the present paper.

\begin{theorem}  \label{th:0.1}  (\cite{Koppelberg:2}, 0.3)  A Boolean
algebra is Cohen 
iff it is the union of a continuous chain  $(A_\alpha)_{\alpha  < \rho}$   
where   $\rho$  is any ordinal, $\pi(A_0) \leq  \omega$,    $A_\alpha$ is a 
regular
subalgebra of 
$A_{\alpha+1}$ and  $\pi(A_{\alpha+1}/A_\alpha ) \leq \omega$.
\end{theorem}
 
The following is a reformulation of a result due to Bandlow (\cite{Bandlow:1}).

\begin{theorem}  \label{th:0.2}    
A Boolean algebra  $A$  is Cohen 
iff  there is a club subset  {\bb S} of  $[A]^\omega$  such that the 
elements
of  {\bb S}  are subalgebras of  $A$  and, for every subset   $T$  of {\bb S}, 
the subalgebra  of  $A$  generated by  $\bigcup T$  is regular in  $A$. 
\end{theorem}

We now give a survey of the proof of our theorem and explain the 
organization of the paper. In fact, what we show is a result on forcing:
we find a Boolean extension  $V^{Q^0}$  of  the universe  $V$  of set
theory which is not Cohen, but some Boolean  extension
$V^{Q^0*Q^1}$  of  $V^{Q^0}$  is. --- After reviewing some material on
Boolean algebras and forcing in Section  1, we
will define forcings  $Q^0$, $Q^1$,  $P^0$, and $P^1$, most of which depend 
on the cardinal  $\kappa$ given in the Main Theorem as a parameter. More 
precisely, we define   $Q^0$  in Section 2 and list some of its basic
properties. In  Section 3, we prove that for  $\kappa \geq \omega_2$,  
$Q^0$   respectively its associated complete Boolean algebra  $B(Q^0)$ 
is not Cohen. We define  $Q^1$ in Section 4,  $P^0$  and   $P^1$  in
Section  5; moreover, we find dense  subsets  $D_Q$  of  the iteration 
$Q^0*Q^1$ of  $Q^0$  and  $Q^1$, respectively  $D_P$  of  $P^0*P^1$, and
prove that  $P^0*P^1$  is Cohen. Finally in Section  6, we 
prove that  $D_Q$  and  $D_P$  are isomorphic. 

This proves the Theorem, 
because of the following well-known facts on the connection between 
partial orderings  $P$  and their associated Boolean algebras  $B(P)$.
For  $D$ a dense subset of  $P$,  $B(D)$  is isomorphic to  $B(P)$;
thus  $B(Q^0*Q^1)$  is  isomorphic to  $B(P^0*P^1)$  and  $B(Q^0*Q^1)$ 
is Cohen. $Q^0$  is completely contained  in the iteration  $Q^0*Q^1$ 
and thus  $B(Q^0)$  is completely embeddable into the  Cohen 
algebra  $B(Q^0*Q^1)$, but   $B(Q^0)$  was not Cohen.

Both iterations  $Q^0*Q^1$  and  $P^0*P^1$  will adjoin the same 
generic objects ($\underline f$  and, for each  $\alpha \in \kappa$,
functions  $\underline t_\alpha  : \omega \rightarrow \omega$  and  
$\underline x_\alpha : \omega \rightarrow 2$),
but in different order; this is why  $B(P^0*P^1)$  is isomorphic to  
$B(Q^0*Q^1)$.
The functions  $\underline t_\alpha$, 
 $\alpha \in \kappa$,  will be almost disjoint
 in the sense that  for $\alpha \neq \beta$,   
 $\underline t_\alpha(i) \neq \underline t_\beta(i)$  will hold for 
almost all  $i$. 
And the generic objects will be connected as follows.
Denote by  $B$  the binary tree of height  $\omega$  and by  $\mbox{lev}_i 
B$ its i\rq th level. For   $i \in \omega$, let  $a_i$  be a subset of 
$\omega$ of  size  $\vert \mbox{lev}_i B \vert = 2^i$ (for ease of
notation,  $a_i$   will later be the set  $\{0, \dots, 2^n -1  \}$, but
this is neither important  nor
necessary for the proof). 
$\underline f$  will be a function from  $B$  into  $\omega$  mapping   
$\mbox{lev}_i B$  onto  $a_i$; for every  $\alpha \in \kappa$,  
$\underline t_\alpha $  and  
$\underline x_\alpha$  will be connected by  $\underline f$  in such a way 
that for almost all  $i$,  $\underline t_\alpha (i) =   \underline f
(\underline x_\alpha \restriction i)$. Now the branches  $\{\underline 
x_\alpha \restriction i :i \in \omega\}$  through  $B$ induced
by the  $\underline x_\alpha$  are almost disjoint, and this causes the 
$\underline t_\alpha$  to be almost disjoint. In fact, this line of argument
reflects the standard
construction of a family of  $2^\omega$  almost disjoint subsets of  
$\omega$ out of the branches of a binary tree:  if we choose the sets 
$a_i$  to be pairwise disjoint, then even the sets 
$\mbox{ran } \underline t_\alpha$,  $\alpha \in \kappa$,  will be almost 
disjoint.

The fact that  $Q^0$  is not Cohen can be partially explained by a 
combinatorial 
principle forced by  $Q^0$, as observed by Soukup (\cite{Juhasz+:1}). Call an 
almost disjoint family  $A$  of subsets of  $\omega$  a  $\kappa$-Luzin
gap if $\vert A \vert = \kappa$  and there is no   $X \subseteq \omega$ 
such  that both $\{a \in A:a \mbox { is almost contained in } X\}$ 
and  $\{a \in A:a \mbox { is almost contained in } \omega \setminus X\}$ 
have size  $\kappa$. Now  $Q^0$  forces that 
$ \{ \mbox{ran } \underline t_\alpha : \alpha \in \kappa \}$  is a 
$\kappa$-Luzin gap. On the other hand, if  $\kappa = \omega_2$  and  CH 
holds in the ground model,  then Cohen\rq s partial order  $Fn(\kappa ,2
)$  forces that  there is no $\kappa$-Luzin gap.

The paper uses basic notions and results on forcing respectively two-step
iterated forcing in the proof of Proposition \ref{pr:3.2} and in Sections
4 and 5. It is however  possible to give an elementary proof which deals
only with partial  orderings and 
their associated Boolean algebras: Proposition \ref{pr:3.2} is provable
without forcing, as  explained
in Section 3. In Section 4, one can define the partial ordering 
$D_Q$  by Proposition \ref{pr:5.5} (plus Proposition \ref{pr:5.4}) and
then check that the map 
$e: Q^0 \rightarrow B(D_Q)$  given by  
$e(p) = \sum^{B(D_Q)} \{(p^\prime, q^\prime) \in D_Q :
p^\prime \mbox{  extends }  p, \mbox{ in }  Q^0\}$ 
is a complete embedding. A similar argument gives a complete embedding  
from  $P^0$  into
$B(D_P)$, and it can be checked in an elementary way that  $B(D_P)$  is 
Cohen.

The first author wants to thank several colleagues for gentle pressure
and constant encouragement during the untimely long preparation of the 
paper, in particular Saka\'e Fuchino, Lutz Heindorf, and Bohuslav
Balcar, and Lajos Soukup for several enlightening remarks on the final
version.

\section{Preliminaries}

For unknown results or unexplained notions, cf.\ \cite{Jech:1}  and  
\cite{Jech:2} 
in set theory, \cite{Koppelberg:1} in Boolean algebras.

\begin{definition}[Boolean algebras]  \label{de:1.1}  The finitary 
Boolean operations are denoted
by  $+$,  $\cdot$,  and  $-$, the infinitary ones by  $\sum$  and  $\prod$. 
$0$  and  $1$  are the distiguished elements. 

For a Boolean algebra  $D$,  $D^+$  is the set  $D \setminus \{0\}$  
of non-zero elements of  $D$. For  $X \subseteq D$,  $\langle X \rangle$  
(respectively  ${\langle X \rangle}^{cm}$)  is the subalgebra of  $D$  
generated by  $X$  (respectively completely generated by  $X$, if  $D$ 
is complete).

$B \leq D$  denotes that  $B$  is a subalgebra of  $D$.  $B$  is a regular 
subalgebra of  $D$  if all infinite sums and products of subsets of  $B$  
that happen to exist in  $B$  are preserved  in  $D$.
\end{definition}

\begin{definition}[Dense subsets of Boolean algebras] \label{de:1.2}  A
subset  $P$  of  $D$ 
is dense in  $D$  if for every  $d \in D^+$  there is  $p \in P$  such that 
$0 < p \leq d$, i.e.\ every element of  $D$  is the least upper bound of some
subset of  $P$. $\pi(D)$, the  $\pi$-weight of  $D$, is the minimal size
of a dense subset of  $D$.  More generally for  $B \leq D$,  the relative  
$\pi$-weight of  $D$  over  $B$, $\pi(D/B)$, is defined as  
$\mbox{min } \{ \vert A \vert : A \subseteq D  \mbox{ and } A \cup B
\mbox{ generates a dense subalgebra of }  D \}$; replacing  $A$  by a 
subalgebra of  $D$  including  $A$, we can assume that  $A$  is a
subalgebra  of  $D$.  In this case  $\{a \cdot b: a \in A, b \in B\}
\setminus \{0\}$  is dense in  $\langle A \cup B \rangle$, so  $\langle
A\cup  B \rangle$   is dense in   $D$  iff, for every  $d \in D^+$,
there are  $a \in A$   and  $b \in B$  such that  $0 < a \cdot b \leq
d$. Moreover, if  $A$  and  $B$ are regular subalgebras of  $D$  and 
$P_A \subseteq A$, $P_B \subseteq B$, $P_D \subseteq D$  are dense  and 
$P_A \cup P_B \subseteq P_D$, then  $\langle A\cup  B \rangle$ is dense 
in  $D$ iff  for every  $q \in P_D$  there are elements  $p \in P_A$  and
$p^\prime \in P_B$  such that  $p \cdot p^\prime > 0$  and for every
$r \in P_D$, if  $r \leq p$  and  $r \leq p^\prime$, then  $r \leq q$;
this is because every element of  $A$  respectively  $B$ is the sum of a
subset of  $P_A$  respectively  $P_B$.

If  $A$  and  $B$  are complete and regular subalgebras of $D$  and  $A
\leq B \leq D$, then  $\pi(B/A) \leq \pi(D/A)$. This is proved as
follows. Choose a subalgebra  $E$  of  $D$  such that  $\vert E \vert  =
\pi(D/A)$  and  $\langle A \cup E \rangle$  is dense in  $D$. Then consider 
the set  $F = \{h(e) : e \in E\}$  where  $h: D \rightarrow B$  denotes
the  projection map given by $h(d) = \mbox{min } \{b \in B : b \geq d\}$ 
from  $D$  to  $B$. Then  $\vert F \vert \leq \pi(D/A)$  and it is easily 
checked
that   $\langle A \cup F \rangle$  is dense in $B$.

For a partially ordered set  $(P, \leq_P)$, we write  $B(P)$  for its 
associated
Boolean algebra or completion, i.e.\  $B(P)$  is the unique complete
Boolean algebra  $B$  such that there is an embedding  $i: P \rightarrow B$
with  $i[P]$   dense in  $B$; cf.\ \cite{Kunen:1} II.3.3. Caused by the
notation on 
forcing 
used in Definition \ref{de:1.3}, we assume that  $i$  is order-reversing.
Moreover,  $i$  is one-
one
and satisfies  $p \leq_P q$  iff  $i(q) \leq_B i(p)$, for all  $p, q \in P$,
iff   $P$  is separative, i.e.\ for  $p$  and  $q$  satisfying  $p \nleq_P q$,
there is  $r \in P$  such that  $q \leq_P r$  and  $r$  is incompatible with  
$p$.
In this case, we will think about  $P$  as being a dense subset of  $B(P)$.
\end{definition}

\begin{definition}[Forcing] \label{de:1.3} 
When dealing with notions of forcing  $(P, \leq)$,  
$p \leq q$  means that the condition  $q$  is stronger than  $p$. For an 
arbitrary
cardinal  $\kappa$,  $Fn(\kappa ,2)$  is the forcing which adjoins  
$\kappa$  Cohen 
reals, i.e.\ a condition in  $Fn(\kappa ,2)$  is a function  $p$  from some 
finite 
subset of  $\kappa$  into  $2$, and  $p \leq q$  holds iff   $p \subseteq q$. 
We 
call its completion  $\mbox{\bb C}_\kappa = B(Fn(\kappa ,2))$  the standard 
Cohen algebra
of  $\pi$-weight  $\kappa$, since  $\pi(\mbox{\bb C}_\kappa) = \kappa$. 
$\mbox{\bb C}_\kappa$ 
is also the completion of the free Boolean algebra over  $\kappa$  
generators; this
is the definition used in  \cite{Koppelberg:2}. As usual in the literature,
we will also 
call
a forcing  $P$  Cohen if  $B(P)$  is isomorphic to some standard Cohen 
algebra.
\end{definition}

\begin{definition}[Two-step iterated forcing]  \label{de:1.4} 
Let us recall the general version 
of 
two-step iterated forcing. If  $(P, \leq_P)$  is a partial ordering in the 
ground
model  $V$  and  $(Q, \leq_Q) \in V^P$  is a  $P$-name for a partial 
ordering 
(i.e.\ if  $P\Vdash ``(Q, \leq_Q) \mbox{ is a partial ordering }"$), 
then the iteration  $P*Q$  of  $P$  and  $Q$  is the partial ordering (in  $V$) 
defined as follows. The elements of  $P*Q$  are certain pairs  $(p,q)$  where  
$p \in P$, $q \in \mbox{dom } Q$  (thus  $q \in V^P$) and  $p \Vdash q \in 
Q$. And  
$(p,q) \leq (p^\prime, q^\prime)$  holds in  $P*Q$  if  $p \leq_P p^\prime$  
and 
$p^\prime \Vdash q \leq_Q q^\prime$. 

When applying this in Sections 4 and 5, we will deal with a simpler 
situation: 
we will have a set  $N$  in  $V$  such that  $P\Vdash Q \subseteq \check 
N$. If
$(p,q) \in P*Q$, let us say that   $p$  decides  $q$  if, for some  $n \in N$,  
$p \Vdash q = \check n$. Here  $\check n$   is the canonical name for  
$n \in V$  in  $V^P$; we will usually write  $n$  for  $\check n$. We call
the subset
\[ \mbox{stp}_N (P*Q) = \{(p,q): p \in P, q \in N, p \Vdash \check q \in Q\} 
\]
of  $P*Q$  the standard part of  $P*Q$  relative to  $N$. This is a dense 
subset
of  $P*Q$, hence their associated Boolean algebras are isomorphic. We will
omit the subscript  $N$  and, if convenient, tacitly pass to a dense subset of
$ \mbox{stp } (P*Q)$  and still call it the standard part of  $P*Q$.
 
Note that if  $p$  decides  $q ( = \check n)$  and if, e.g.,   
$n$  is, in  $V$, a function, then  $p$  also decides the domain, the range, 
and the values of  $q$  since the statements  $``u = \mbox{dom } n"$,  
$``v = \mbox{ran } n"$,  $``i \in u = \mbox{dom } n \mbox{ and } j = n(i)"$  are 
$\Delta_0$, hence absolute for  $V$  and  $V^P$.
\end{definition}
 
\begin{definition}[Some definitions for Sections 2 to 6] \label{de:1.5}   
We fix some
notation which will be used throughout the paper. 

For $n \in \omega$, let  $a_n = 2^n = \{0, \dots, 2^n - 1\} \subseteq 
\omega$. 

$T$  is the tree  of height  $\omega$  with  $n$\rq th  level
$\mbox{lev}_n T$  the set of those functions
$t$  from  $n$  to  $\omega$  such that, for  $i < n$,  $t(i) \in a_i$, i.e.
$\mbox{lev}_n T = a_0 \times \dots \times a_{n-1}$. 
$B$  is the binary tree of height  $\omega$  with $\mbox{lev}_n B$  
the set of all functions from  $n$  to  $2$. In both cases, the tree ordering is
set-theoretic inclusion. 

For  $M$  a set of sequences with common domain some  $n \leq \omega$  
and   $k \leq n$, we say that the elements of  $M$  are disjoint  
 above  $k$  if for every  $i \in [k,n)$,
the values  $m(i)$,  $m \in M$, are pairwise distinct.
\end{definition}

\section{Simple properties of $Q_X^0$}

Our investigation of the forcing  $Q^0$  will use, more generally, 
the forcings  $Q_X^0$  where  $X$  is an arbitrary set (or, in section 3, 
a subset of some cardinal  $\kappa \geq \omega_2$). 
$Q^0$  will simply be the special case  $Q_X^0$  
where  $X = \kappa$. In this section, we collect some basic properties of  
the
$Q_X^0$.

\begin{definition}  \label{de:2.1}
 For any set  $X$, the forcing  $Q_X^0$  is defined 
as 
follows. An element of  $Q_X^0$  is  a function  $p$  with  $\mbox{dom } p$  
a finite subset  of  $X$   such that, 
for some  $n \in \omega$ (the {\it height} of   $p$,  $\mbox{ht } p$):

(a) $\vert \mbox{dom } p  \vert \leq a_n$

(b) writing $t_\alpha^p$  for  $p(\alpha)$: each  $t_\alpha^p$, 
$\alpha \in \mbox{dom } p$, is an element of  $\mbox{lev}_n T$,
i.e.\  $t_\alpha^p(i) \in a_i$  for  $i < n$  

(c) $\mbox{dom } p = \emptyset$  
implies   $\mbox{ht }p  = 0$ and  $\mbox{dom } p \neq \emptyset$  
implies   $\mbox{ht }p  \geq 2$  (this is only for technical convenience).
 
For  $p$  and  $q$  in  $Q_X^0$, $p \leq q$  iff  

(d) $\mbox{dom } p \subseteq \mbox{dom } q$  and  
$\mbox{ht } p \leq \mbox{ht } q$

(e)  for  $\alpha \in \mbox{dom } p$,  $t_\alpha^p \subseteq t_\alpha^q$

(f)  the  $t_\alpha^q$,  $\alpha \in \mbox{dom } p$, are  disjoint
above $\mbox{ht } p$.
\end{definition}

\begin{lemma}[and Definition] \label{le:2.2}
For every  $k \in \omega$, the set  
$\{q : \mbox{ht } q \geq k \}$  is dense in  $Q_X^0$. Also for every  
$\alpha \in X$,   $\{q : \alpha \in \mbox{dom } q\}$  is dense in  $Q_X^0$.
It follows that, if  $G \subseteq Q_X^0$  is  $Q_X^0$-generic over the ground
model  $V$ and we put
$ {\underline t}_{\alpha G} = \bigcup \{t_\alpha^p : p \in G
\mbox{ and } \alpha \in \mbox{dom } p\}  $
for  $\alpha \in X$, each  ${\underline t}_{\alpha G}$  is an element of
the cartesian product  $\prod_{i \in \omega} a_i$.
\end{lemma}

\proof   Obvious.  For the first claim, note that 
this is where (a) of Definition \ref{de:2.1} is
used; for the second one, enlarge first the height of a given condition in  
$Q_X^0$, if necessary, and then the domain.    \endproof

The subsequent propositions will use the following criterion 
for compatibility in  $Q_X^0$.

\begin{proposition} \label{pr:2.3}  Assume  $p$  and  $q$  are in  $Q_X^0$ 
and  $\mbox{ht } p \leq \mbox{ht } q$. 
Then  $p$  and  $q$  are compatible in  $Q_X^0$  iff

(a)  for   $\alpha \in \mbox{dom } p \cap \mbox{dom } q$, 
$t_\alpha^p \subseteq   t_\alpha^q$

(b)  the  $t_\alpha^q$,  $\alpha \in \mbox{dom } p \cap \mbox{dom } q$,
are  disjoint above  $\mbox{ht } p$.
\end{proposition}

\proof  Obvious.     \endproof

\begin{corollary} \label{co:2.4}  
$Q_X^0$  satifies the countable chain condition.           
\end{corollary}

\proof  By the usual $\Delta$-system lemma argument.      \endproof

The next proposition says that we can think of  $Q_X^0$  as  being 
a dense subset of its associated Boolean algebra  $B_X = B(Q_X^0)$.

\begin{proposition}  \label{pr:2.5} 
$Q_X^0$  is separative.            \end{proposition}

\proof  Assume  $p$  and  $q$  are in  $Q_X^0$  and 
$p \nleq q$. We will find  $r \in Q_X^0$  such that  $q \leq r$  and
$r$  is incompatible with  $p$. Without loss of generality, 
we may assume that  $p$  and  $q$  are compatible.           

{\it Case 1.}   $\mbox{dom } p \nsubseteq \mbox{dom } q$. Fix  $\alpha \in
\mbox{dom } p \setminus \mbox{dom } q$. 
We will choose  $r \geq q$  such that  $\mbox{dom } r = \mbox{dom } q
\cup \{\alpha \}$,  $\mbox{ht } r > \mbox{ht } q$,  and  
$\vert \mbox{dom } r \vert \leq a_{\mbox{ht } r}$. More precisely, take  
$t_\alpha^r \in a_0 \times \dots \times a_{\mbox{ht } r -1}$  
such that  $t_\alpha^r(1) \neq t_\alpha^p(1)$  (recall Definition
\ref{de:2.1}(c) and
$a_1 = 2$). For  $\beta \in \mbox{dom } q$, let  
$t_\beta^r \in \mbox{lev}_{\mbox{ht } r} $  extend  $t_\beta^q$  in 
such a way that all  $t_\beta^r(i)$, $\mbox{ht } \leq i <  \mbox{ht } r$, 
are distinct.

{\it Case 2.}  $\mbox{dom } p \subseteq \mbox{dom } q$, but  
$\mbox{ht } q < \mbox{ht } p$. $\mbox{ht } p > 0$, so by Definition 
\ref{de:2.1}(c), also  $\mbox{dom } p \neq \emptyset$; fix  $\alpha \in
\mbox{dom } p$. Now take 
$r \geq q$  such that  $\mbox{dom } r =  \mbox{dom } q$, 
$\mbox{ht } r = \mbox{ht } p$, and  $t_\alpha^r(\mbox{ht } q) \neq 
t_\alpha^p(\mbox{ht } q)$; thus  $r$  and  $p$  are incompatible.
This is possible since  $\mbox{dom } q \neq \emptyset$  
(so  $\mbox{ht }q \geq 2 $)  and  $a_{\mbox{ht } q} \geq 2$.

{\it Case 3.}  $\mbox{dom } p \subseteq \mbox{dom } q$  and  
$\mbox{ht } p \leq \mbox{ht } q$.  $p$  and  $q$  are compatible, 
so by Proposition \ref{pr:2.3}, we have that for every  $\alpha \in
\mbox{dom } p $,
$t_\alpha^p \subseteq t_\alpha^q $  and the  $t_\alpha^q$,  
$\alpha \in \mbox{dom } p$, are  disjoint above $\mbox{ht } p$.
But then  $p \leq q$, a contradiction.
\endproof

In the two subsequent propositions, we use the following construction. 
Every
permutation  $h$  of  $X$  induces an automorphism  $\overline h$  of the
partial ordering  $Q_X^0$  by letting, for  $p \in Q_X^0$,
$ \mbox{dom } \overline h p  = h [\mbox{dom } p]$,
$ \mbox{ht } \overline h p  = \mbox{ht } p  $, and, 
for  $\alpha \in  \mbox{dom } p$,
$ t_{h(\alpha)}^{ \overline h p} =  t_\alpha^p $.
Moreover, we call a forcing  $Q$  {\it weakly homogeneous} if for 
arbitrary  $p$  and  $q$  in  $Q$, there are  $p^{\prime}$  and  
$q^{\prime}$  in  $Q$  such that  $p \leq p^{\prime}$,  $q \leq q^{\prime}$,
and  $Q \restriction p^{\prime}$  is isomorphic to  $Q \restriction 
q^{\prime}$,
where  $Q \restriction r =  \{x \in Q : x \geq r\}$, for  $r \in Q$.

\begin{proposition} \label{pr:2.6}   $Q_X^0$  is weakly homogeneous, for
infinite  
$X$.
\end{proposition}

\proof  Let  $p$  and  $q$  in  $Q_X^0$  be given. Fix a permutation
$h$  of  $X$  such that  $h [\mbox{dom } q]$  is disjoint from 
$\mbox{dom } p$. By Proposition \ref{pr:2.3}, there is a common extension 
$r$ of  $p$  and  
$\overline h (q)$. Now  $r \geq \overline h (q)$,  
$\overline {h^{-1}} (r) \geq q$,  and  $Q_X^0 \restriction r$  is 
  isomorphic to   $Q_X^0 \restriction \overline {h^{-1}} (r)  $.  \endproof

\begin{proposition}  \label{pr:2.7}   
Assume $X$  is infinite and  $X\subseteq Y$. 
Then
$Q_X^0$  is completely contained in  $Q_Y^0$  (cf.\ \cite{Kunen:1} VII.7.1  
for this notion).
\end{proposition}

\proof   Clearly,  $Q_X^0$  is a subordering of  
$Q_Y^0$, and, by Proposition \ref{pr:2.3}, two elements of  
$Q_X^0$  are compatible in  $Q_X^0$  iff they are in  $Q_Y^0$.Thus 
assume  $p^\prime \in Q_Y^0$  with the aim of finding 
$p \in Q_X^0$  such that every extension of  $p$  in  $Q_X^0$  is 
compatible wih  $p^{\prime}$.

Write  $\mbox{dom } p^{\prime} = r \cup s^{\prime}$  where  
$r \subseteq X$  and  $s^{\prime} \subseteq Y \setminus X$.
Then choose a subset $s$  of  $X$  disjoint from  $r$  such that 
$\vert s \vert = \vert s^{\prime} \vert$  and a 
permutation  $h$  of  $Y$  satisfying  $h \restriction r = id$  and
$h[s^{\prime}] = s$. We will show that  $p = \overline h (p^{\prime})$  
works for our claim; note that  $\mbox{dom } p = r \cup s$.

In fact, assume that  $q \in Q_X^0$  extends  $p$; say 
$\mbox{dom } q = r \cup s \cup u$  where  $u$  is disjoint from  $r \cup s$.
Choose another permutation  $k$  of  $Y$  such that  $k$  and  $h^{-1}$
coincide on  $\mbox{dom } p$  and  $k$  maps  $u$  onto a subset  
$u^{\prime}$ 
of  $Y$  disjoint from  $\mbox{dom } q$. Clearly  $q^{\prime} = \overline k 
(q)$ 
extends  $p^{\prime}$; thus it suffices to prove that  
$q$  and  $q^{\prime}$  are compatible. But  $\mbox{ht } q 
= \mbox{ht } q^{\prime}$, $\mbox{dom } q  \cap  \mbox{dom } q^\prime = r$,
and for every  $\alpha \in r$, we have  $k(\alpha) = \alpha$  and thus 
 $t_\alpha^q = t_\alpha^{q^\prime}$.     \endproof

\section{$Q^0$  is not Cohen}

We prove in this section that the forcing  $Q^0 = Q_\kappa^0$  is not Cohen, 
for  $\kappa \geq \omega_2$, i.e.\ its associated Boolean algebra  $B(Q^0)$   
is not Cohen.
The ideas lying behind the proof are from \cite{Koppelberg:2}  (cf.\
Theorem 0.1 in 
the introduction), i.e.\ essentially from Shapiro\rq s proof that subalgebras 
of free
algebras are Cohen, but we give a completely self-contained presentation 
here.
The main argument in the proof is the following lemma. We have not tried 
to minimize its assumptions since they are so naturally satisfied 
in the intended application.

\begin{lemma} \label{le:3.1}  Assume that  $\kappa \geq \omega_2$  is a 
cardinal
and that, for every subset  $X$  of  $\kappa$,
we are given two separative partial orderings  $P_X$  and  $Q_X$  with the 
following properties. We write  $A_X = B(P_X)$,  $B_X = B(Q_X)$, and 
assume
that  $P_X$  is a dense subset of  $A_X$; similarly for  $Q_X$  and  $B_X$. 

(a) $P_X$  and  $Q_X$  satisfy the countable chain condition

(b)  $X \subseteq Y \subseteq \kappa$  implies that  $P_X \subseteq P_Y$  
and  $P_X$  is completely contained  in 
$P_Y$  (so without loss of generality,  $A_X$  is a 
regular subalgebra of  $A_Y$); similarly for  $Q_X \subseteq Q_Y$  and 
$B_X \leq B_Y$

(c)  $P_X = \bigcup \{P_e : e \subseteq X \mbox{ finite}\}$; 
similarly for  $Q_X$

(d) $\vert P_X \vert \leq \vert X \vert$  for infinite  X; similarly for  $Q_X$

(e)  for  $X, Y \subseteq \kappa$,  $A_{X \cup Y}$  is completely generated 
by  $A_X \cup A_Y$;  similarly for  $B_{X \cup Y}$

(f) if  $Y$  is countable, then  $\pi(A_{X \cup Y}/A_X) \leq \omega$.

Assume  $B_\kappa$  is isomorphic to  $A_\kappa$.  Then there is a club 
subset
{\bb C}  of  $[\kappa]^{\omega_1}$  such that
 
 (g)  for  $X \in \mbox{\bb C}$  and  $Y \subseteq \kappa$  countable,
 $\pi(B_{X \cup Y}/B_X) \leq \omega$.

\end{lemma}

\bigskip

We will apply Lemma \ref{le:3.1} to the situation where  $P_X = Fn(X,2)$ 
is standard 
Cohen 
forcing and  $Q_X = Q_X^0$  as defined in Definition \ref{de:2.1}. 
(a) through (f) of Lemma \ref{le:3.1} are 
clearly
satisfied for the forcings  $P_X$, and (a) through (d) hold for  $Q_X$,
by the results of Section 2. We prove in the subsequent lemmas that the 
$Q_X$
satisfy (e), but not (g) --- hence  $B(Q_\kappa^0)$  is not isomorphic to 
$B(Fn(\kappa,2))$.

\proofof of Lemma \ref{le:3.1}.  For convenience of notation, we assume
that 
$A_\kappa$  and  $B_\kappa$  are the same Boolean algebra  $D$;  so for 
$X \subseteq \kappa$, $A_X$  and  $B_X$  are regular subalgebras of  $D$.

Call a subset  $M$  of  $D$  nice if there is a regular complete subalgebra
$C$  of  $D$  such that  $M$  is  a dense subset of  $C$. (E. g.,  $P_X$  and  
$Q_X$  are nice.) In this case,   $C$  is uniquely determined  by  $M$  since 
$C = \{\sum^D M_0: M_0 \subseteq M\}  $
and we write  $C = C(M)$.

For  $M$  and  $N$  subsets of  $D$, we say that  $M$  is dense for  $N$
if for every  $n \in N$, there is  $M_0 \subseteq M$  such that  $n = \sum^D 
M_0$.
Thus if  $M$  and  $N$  are nice,  $C(M) = C(N)$  iff  $M$  is dense for 
$N$  and  $N$  is dense for  $M$.

Define 
\[ \mbox{\bb C} = \{X \subseteq \kappa: \vert X \vert = \omega_1, P_X \mbox{ is 
dense for } Q_X \mbox{ and } Q_X \mbox{ is dense for } P_X\}. \]
It follows from the assumptions (a), (c), and (d) that  
{\bb C}  is club in  $[\kappa]^{\omega_1}$. And for  $X \in \mbox{\bb C}$, we 
have 
$A_X = B_X$.

To check (g), assume  $X$  is an element of {\bb C} and  
$Y \subseteq \kappa$  is countable. By (a), (c), and (d) again, we find 
a countable  $Z \subseteq \kappa$  such that  $P_Z$  is dense for  $Q_Y$; so  
$B_Y \leq A_Z$. Now by (e),

\[ A_X = B_X \leq B_{X \cup Y} =
{\langle B_X \cup B_Y \rangle}^{\mbox{ cm}}
\leq {\langle A_X \cup A_Z \rangle}^{\mbox{ cm}}   
= A_{X \cup Z}  \]
and  $\pi(A_{X \cup Z}/A_X) \leq \omega$  holds by (f). It follows from 1.2  
that
\[ \pi(B_{X \cup Y}/B_X) \leq \omega.  \]                 \endproof

For the rest of this section, fix a cardinal  $\kappa \geq \omega_2$  
and write, for  $X \subseteq \kappa$:
\[ B_X = B(Q_X^0),  \]
a regular subalgebra of  $B_\kappa$. We proceed to show that assumption 
(e) 
of Lemma \ref{le:3.1} holds for the  $B_X$, but (g) fails. In  $B_\kappa$,
note that  $B_X$  
has  $Q_X^0$  as set of complete generators, since  $Q_X^0$  is a dense 
subset
of  $B_X$  and thus every element of  $B_X$  is  (in $B_\kappa$) the join 
of some subset of  $Q_X^0$.

The details of the following proof may be unnecessary for a reader 
experienced
with forcing. On the other hand, the use of forcing can be avoided, at the
price of a little more computation and less insight: simply {\it define} 
$b_{\alpha i j}$
to be  $\sum N_{\alpha i j}$  (where  $N_{\alpha i j}$  is as in the proof
of Claim 1 below); 
this makes Claim 1 trivial. For Claim 2, just prove that, for  $p \in Q_X^0$ 
with domain  $u$  and height  $n$:
\[ p = \prod \{b_{\alpha i j}: \alpha \in u, i < n, t_\alpha^p(i) = j  \} 
\cdot
\prod \{-(b_{\beta i j} \cdot b_{\gamma i j}): \beta \neq \gamma \mbox{ 
in } u, 
i \geq n, j \in a_i\}. \]

\begin{proposition}  \label{pr:3.2} 
For  $X, Y \subseteq \kappa$,  
$B_{X \cup Y}$  is completely generated by  $B_X \cup B_Y$.            
\end{proposition}

\proof   For  $\alpha \in \kappa$,  $i \in \omega$,  and  $j \in a_i$, 
let  
$\sigma_{\alpha i j}$  be the sentence  ``${\underline t}_\alpha(i) = j$" 
of the forcing language over  $Q_\kappa^0$ 
and let  $b_{\alpha i j}$  be its Boolean truth value  
$ \Vert \sigma_{\alpha i j} \Vert $, computed in  $B_\kappa$.
Here  ${\underline t}_\alpha$  is (the canonical name for)
the generic object introduced in Lemma \ref{le:2.2}; we do not distinguish
notationally 
between 
$\alpha, i, j$  in the ground model and their canonical names 
$\check \alpha, \dots $  in the forcing language.

For  $X \subseteq \kappa$, consider the set  
$M_X = \{ b_{\alpha i j}: \alpha \in X,  i \in \omega, j \in a_i\}$, 
and the subsequent claims.

{\it Claim 1.}  $M_X \subseteq B_X$.

{\it Claim 2.}  $M_X$ completely generates  $B_X$.

Now clearly  $M_{X \cup Y} = M_X \cup M_Y$, and thus the proposition
follows from the claims.

To prove the claims, recall the basic fact on the connection between forcing
and Boolean-valued models (cf.\ \cite{Jech:1} p.\ 166): for  $p \in
Q_\kappa^0$  
and 
$\sigma$  a sentence of the forcing language over  $Q_\kappa^0$, 
\[ p \Vdash  \sigma \mbox{ iff } p \leq_B  \Vert  \sigma \Vert   .\]

\newpage

Here we write  $\leq_B$  for the Boolean partial ordering. All joins and 
meets
 below are computed in  $B_\kappa$.

{\it Proof of Claim 1.}  We prove that for  $\alpha \in X$, $i \in \omega $ , 
$j \in a_i$,
we have $b_{\alpha i j} = \sum N_{\alpha i j}$  where  
$N_{\alpha i j} = \{p \in Q_{\{\alpha\}}^0 :
\mbox{dom } p = \{\alpha\},\mbox{ht } p > i, t_\alpha^p(i) = j \}$.
Here,  $\geq$  is obvious since every  $p \in  N_{\alpha i j}$  forces  
$\sigma_{\alpha i j}$. 
Assume for contradiction that  $b_{\alpha i j}$  is (in  $B_\kappa$) 
strictly greater than  $\sum N_{\alpha i j}$. Then there is some  
$q \in Q_\kappa^0$  
forcing  $\sigma_{\alpha i j}$ but incompatible with all  $p \in N_{\alpha i 
j}$.
By extending  $q$, we may assume that  $\alpha \in \mbox{dom } q$ 
and  $i < \mbox{ht } q$.
Consider $k = t_\alpha^q(i)$. 
Now  $k \neq j$: otherwise, let  $p$  be the restriction of  $q$  with domain  
$\{ \alpha\}$  and  height  $i+1$.        
Then   $p \in N_{\alpha i j}$  and  $q$  extends  $p$  in  $Q_\kappa^0$;  
a contradiction since  $q$  was incompatible with all  $p \in N_{\alpha i j}$. 
-
It follows that  $q \Vdash {\underline t}_\alpha(i) = k \neq j$  and 
thus  $q \Vdash \neg \sigma_{\alpha i j}$, a contradiction.

{\it Proof of Claim 2.}  It suffices to show that  $Q_X^0 \subseteq B_X$  is 
completely generated by  $M_X$. So let  $p \in Q_X^0$, say with domain  
$u$  and height  $n$, and consider the set of sentences of the forcing 
language 
\[ \Sigma_p = \{ \sigma_{\alpha i j}: \alpha \in u, i < n, t_\alpha^p(i) = j  
\} 
\cup
\{ \neg(\sigma_{\beta i j} \land \sigma_{\gamma i j}): 
                        \beta \neq \gamma \mbox{ in } u, i \geq n, j \in
a_i\}.\]

The Boolean value of each  $\sigma \in \Sigma_p$  is clearly generated by  
$M_u \subseteq M_X$, thus it suffices to prove that  
$p = \prod \{ \| \sigma \| : \sigma \in \Sigma_p \}$, i. e. that for each  
$q \in Q_\kappa^0$,  $q$  extends  $p$  iff  $q \Vdash \Sigma_p$ 
(where  $q \Vdash \Sigma_p$  means that  $q \Vdash \sigma$, for every  
$\sigma \in \Sigma_p$). Here,  $\Rightarrow$  is clear since  $p \Vdash 
\Sigma_p$,
 by  (f)  of Definition \ref{de:2.1}. Conversely, assume for contradiction
that  $q \Vdash 
\Sigma_p$  
 but  $q$   does not extend $p$. By applying Proposition \ref{pr:2.5}  and
extending  $q$, we 
may assume
 that  $q$  is incompatible with  $p$,  $ u \subseteq  \mbox{dom } q$  
 and  $n \leq \mbox{ht } q$.
 By Proposition \ref{pr:2.3}, we have to distinguish two cases. Either there are
 $\alpha \in u$  and  $i < n$  such that  
 $t_\alpha^p(i) \neq t_\alpha^q(i)$; then  $q \Vdash \neg \sigma_{\alpha i 
j} $
 where  $j = t_\alpha^p(i)$, 
 contradiction.  Or there are   $i \in [n, \mbox{ht } q)$  and  $\beta \neq 
\gamma$
 in  $u$  such that  $t_\beta^q(i) = t_\gamma^q(i)$; then 
 $q \Vdash \sigma_{\beta i j} \land \sigma_{\gamma i j}$  where  $j = 
t_\beta^q(i)$,
 a contradiction again. 
 \endproof

 The example given in Proposition \ref{pr:3.4} below is the crucial fact
responsible for the 
failure of
 (g) in Lemma \ref{le:3.1} for the algebra  $B_\kappa$. We need another
easy lemma on 
the forcings
 $Q_X^0$  for this. 
 
\begin{lemma}  \label{le:3.3}  
Let  $X$  and  $Y$  be arbitrary sets and assume 
that
  $p \in Q_X^0$  and  $p^\prime \in Q_Y^0$  are compatible in  $Q_{X \cup 
Y}^0$,  
  $k \in \omega$,  
  $\alpha \in X$, and  $\beta \in Y$. Then there are compatible  $q \in 
Q_X^0$  
  and  $q^\prime  \in Q_Y^0$  such that:  $p \leq q$,   $p^\prime \leq 
q^\prime$,
  $\alpha \in  \mbox{dom } q$,  $\beta \in \mbox{dom }q^\prime $, and  
  $\mbox{ht } q  = \mbox{ht } q^\prime \geq k$.
\end{lemma}

\proof  In  $Q_{X \cup Y}^0$, take a common extension  $r$  of  $p$  
and  
$p^\prime$.
By extending  $r$, we may assume that  $\mbox{ht } r \geq k$  and  
$\alpha, \beta \in \mbox{dom } r$. Then let  $q$  respectively  $q^\prime$ 
be the restrictions of  $r$  to  $\mbox{dom } p \cup \{\alpha\}$  
respectively  $\mbox{dom } p^\prime \cup \{\beta\}$. \endproof

\begin{proposition}  \label{pr:3.4} Let  $T$  be a proper subset of  $X$ 
and let  $Y$ 
be a non-empty set disjoint from  $X$. Then, for some  $q \in Q_{X \cup 
Y}^0$, 
there are no compatible  $p \in Q_X^0$  and  $p^\prime \in Q_{T \cup Y}^0$  
such that every common extension of  $p$  and  $p^\prime$  extends  $q$.   
       
\end{proposition}

\proof  Fix elements  $\alpha \in X \setminus T$,  $\beta \in Y$;  
thus   $\alpha \in X$,  $\beta \in Y$, and  $\alpha \neq \beta$. 
Let  $q$  be an arbitrary element of  $Q_{X \cup Y}^0$  satisfying  
$\mbox{dom } q= \{\alpha, \beta \}$.

Assume for contradiction that we have compatible  $p \in Q_X^0$  and  
$p^\prime \in Q_{T \cup Y}^0$  such that every common extension of  $p$  
and
$p^\prime$  extends  $q$. Applying Lemma \ref{le:3.3} and extending  $p$ 
and  
$p^\prime$  if 
necessary, we may assume that  $\alpha \in  \mbox{dom } p$, 
$\beta \in  \mbox{dom } p^\prime$, $\mbox{ht } p = \mbox{ht }p^\prime
 = m \geq \mbox{ht } q$, and  
 $w = \mbox{dom } p \cup \mbox{dom } p^\prime$  has size at most  $a_m$.
 
We choose a common extension  $r \in Q_{X \cup Y}^0$  of  $p$  and  
$p^\prime$  
as follows: put  $\mbox{dom } r = w$  and  $\mbox{ht } r = m+1$. We are left 
with
defining  $t_\gamma^r(m)$, for all  $\gamma \in w$. Simply define these 
values 
such that: all  $t_\gamma^r(m)$,  $\gamma \in \mbox{dom } p$, are 
distinct; all  $t_\gamma^r(m)$,  $\gamma \in \mbox{dom } p^\prime$, are 
distinct;
but  $t_\alpha^r(m) = t_\beta^r(m)$. This is possible since  
$\alpha \in \mbox{dom } p \setminus \mbox{dom } p^\prime$  and
$\beta \in \mbox{dom } p^\prime \setminus \mbox{dom } p$.

By our assumption above,  $r$  must extend  $q$. But this is not the case, 
since
$\alpha, \beta$  are distinct elements of  $\mbox{dom } q$, 
$m \in [\mbox{ht } q, \mbox{ht } r)$, and  $t_\alpha^r(m) = t_\beta^r(m)$. 
        
\endproof

\begin{corollary}  \label{co:3.5} Assume  $X \subseteq \kappa$  is uncountable 
and  
$Y \subseteq \kappa$  is nonempty and disjoint from  $X$.  Then  
$\pi(B_{X \cup Y}/B_X)$  is uncountable.            \end{corollary}

\proof   If not, we can find a countable subset  $C$  of  $B_{X \cup 
Y}$  
such that  $B_X \cup C$  generates a dense subalgebra of  $B_{X \cup Y}$. 
Choose a countable subset  $T$  of  $X$  such that  $C \subseteq B_{T \cup 
Y}$; 
the subalgebra generated by  $B_X \cup B_{T \cup Y}$  is still dense in  
$B_{X \cup Y}$. A remark in 1.2, applied to  $B_X,  B_{T \cup Y}  \leq
B_{X \cup Y}$, gives that for every  $q \in Q_{X \cup Y}^0$, there are 
compatible 
$p \in Q_X^0$  and  $p^\prime \in Q_{T \cup Y}^0$  such that every 
common
extension of  $p$  and  $p^\prime$ (in  $Q_{X \cup Y}^0$) extends  $q$.
But this contradicts Proposition \ref{pr:3.4}.
\endproof

\begin{theorem}  \label{th:3.6}
 $B_\kappa$  ($= B(Q_\kappa^0)$)  is not Cohen,
for  $\kappa \geq \omega_2$.            \end{theorem}

\proof  Assume it is. It is easily checked that  $\pi(B_\kappa) = 
\kappa$.
Also  $B_\kappa$  is weakly homogeneous, by Proposition \ref{pr:2.6}; thus 
$B_\kappa$  must 
be
isomorphic to the standard Cohen algebra  $\mbox{\bb C}_\kappa$  of  $\pi$-
weight  $\kappa$.           

By Lemma \ref{le:3.1}, there is some  $X \in [\kappa]^{\omega_1}$  such
that for every 
countable
$Y \subseteq \kappa$,  $\pi(B_{X \cup Y}/B_X) \leq \omega$, 
contradicting Corollary \ref{co:3.5}. 
\endproof

\section{$Q^1$  and a dense subset of $Q^0 * Q^1$}

For the remaining sections, let  $\kappa$  be  an arbitrary cardinal and put, 
as before,  $Q^0 = Q_\kappa^0$. Following the plan in the introduction, we 
will
define a forcing  $Q^1$  in  $V^{Q^0}$, i.e.\ a  $Q^0$-name  $Q^1$  such that
$Q^0 \Vdash ``Q^1 \mbox{ is a partial ordering }"$, describe the standard part
of  $Q^0 * Q^1$  and find a dense subset of the standard part. The definition 
of 
$Q^1$  will 
use the generic functions  ${\underline t_\alpha}_G : \omega \rightarrow 
\omega$ 
defined in Lemma \ref{le:2.2} for  $\alpha \in \kappa$  respectively their
canonical $Q^0$-
names 
$\underline t_\alpha$.  Recall from Definition \ref{de:1.5} the definitions
concerning the trees  $B$  and  $T$ and the numbers  $a_n$.

\begin{definition}  \label{de:4.1}  (In  $V^{Q^0}$) An element of  $Q^1$ 
is  a pair  
$q = ( f^q ,(x_\alpha^q)_{\alpha \in u}  )$ such that for some finite  
$u \subseteq \kappa$  (the domain of $q$,  $\mbox{dom } q$) 
and some  $m \in \omega$ (the height of   $q$,  $\mbox{ht } q$):

(a)  $f^q$  maps  $\bigcup_{i<m} \mbox{lev}_i B$  into  $\omega$; for  $i<m$, 
$f^q \restriction \mbox{lev}_i B$  is a bijection from  $\mbox{lev}_i B$  onto  
$a_i$  (note  $a_i = 2^i = \vert \mbox{lev}_i B  \vert$)

(b) the  $x_\alpha^q$,  $\alpha \in u = \mbox{dom } q$, are pairwise distinct
elements of  $\mbox{lev}_m B$ (and thus  $\vert \mbox{dom } q  \vert 
\leq 2^m = a_m$) 

(c) the  $\underline t_\alpha$,  $\alpha \in u = \mbox{dom } q$,  are  
disjoint above  $m = \mbox{ht } q$. 

For  $q$  and  $q^\prime$  in  $Q^1$, $q \leq q^\prime$  iff  

(d) $\mbox{dom } q \subseteq \mbox{dom } q^\prime$  and  
$\mbox{ht } q \leq \mbox{ht } q^\prime$

(e)  $f^q \subseteq f^{q^\prime}$  and, for  $\alpha \in \mbox{dom } q$,  
$x_\alpha^q \subseteq x_\alpha^{q^\prime}$

(f)  for  $\alpha \in \mbox{dom } q$  and  $i \in [\mbox{ht } q,
\mbox{ht } q^\prime)$,  ${\underline t}_\alpha(i) = f^{q^\prime}
(x_\alpha^{q^\prime} \restriction i)$.
\end{definition}

For  $H$  $Q^1$-generic over  $V^{Q^0}$, we obtain the generic objects  
${\underline f}_H = \bigcup \{f^q : q \in H\}$, a map from
the binary tree  $B$  into  $\omega$  which maps the  $i$\rq th level of  
$B$ 
in a one-one manner onto  $a_i$, and, for  $\alpha \in \kappa$,
${\underline x}_{\alpha H} = \bigcup \{x_\alpha^q : 
q \in H \mbox{ and } \alpha \in \mbox{dom } q \} : \omega \rightarrow 2$.
They are related to the generic objects  $\underline t_{\alpha G}$  adjoined 
by  $Q^0$  by the fact that, for almost all  $i \in \omega$ 
(we omit the subscripts  $G$  and  $H$),  
$\underline t_\alpha (i) =   \underline f (\underline x_\alpha \restriction 
i)$.

Let us describe the standard part  $\mbox{stp } (Q^0*Q^1)$  of  $Q^0*Q^1$. All 
conditions defining the elements of respectively the partial order on $Q^1$ 
in Definition \ref{de:4.1} deal with objects in  $V$  or 
are absolute, except (c). But for  $p \in Q^0$,  $u \subseteq \mbox{dom } p$ 
finite, 
and  $m \leq \mbox{ht } p$,  
$p$ forces the  $\underline t_\alpha$,  $\alpha \in u$, to be  
disjoint above  $m$  iff the  $t_\alpha^p$,  $\alpha \in u$, are  disjoint 
above  $m$  (i.e.\  disjoint on the interval  $[m,\mbox{ht } p)$). 
This is because  $p$  forces that the  $\underline t_\alpha$,  $\alpha \in 
u$, are
  disjoint above  $\mbox{ht } p$.

\begin{proposition} \label{pr:4.2} The elements of  $\mbox{stp } (Q^0*Q^1)$
 are 
those pairs  
$(p,q)$  such that  $p \in Q^0$   and  $q$  is a pair  
$( f^q ,(x_\alpha^q)_{\alpha \in u})$  where  $u \subseteq \mbox{dom } p$,  
$m = \mbox{ht } q \leq \mbox{ht } p$  satisfying (in  $V$)   (a) and
(b) of Definition \ref{de:4.1}, plus

(c) the  $t_\alpha^p$,  $\alpha \in u$, are  disjoint above  $m$.

For  $(p,q)$  and  $(p^\prime, q^\prime)$  in  $\mbox{stp } (Q^0*Q^1)$,  
$(p,q) \leq (p^\prime, q^\prime)$  iff  
(d) and (e) of Definition \ref{de:4.1} hold, plus

(f) for  $\alpha \in \mbox{dom } q$  and  $i \in [\mbox{ht } q,
\mbox{ht } q^\prime)$,  ${\underline t}_\alpha^{p^\prime}(i) = 
f^{q^\prime}(x_\alpha^{q^\prime} \restriction i)$.
\end{proposition}

\begin{proposition} \label{pr:4.3} The following subset of  $\mbox{stp }
(Q^0*Q^1)$  
is dense in  $\mbox{stp } (Q^0*Q^1)$, hence in  $Q^0*Q^1$.
\[ D_Q = \{(p,q) \in \mbox{stp } (Q^0*Q^1) : \mbox{dom } p = \mbox{dom } q  
\mbox{ and }  \mbox{ht } p = \mbox{ht } q\}. \]
\end{proposition}

\proof  Let  $(p,q) \in \mbox{stp } (Q^0*Q^1)$  be given; we find  
$q^\prime$ 
such that  $(p,q^\prime)$  is in  $D_Q$  and extends  $(p,q)$. Write  
$v = \mbox{dom } p \supseteq u = \mbox{dom } q$,  $n = \mbox{ht } p \geq m
= \mbox{ht } q $.

First pick, for  $\alpha \in v$, an element  $x_\alpha^{q^\prime}$  of 
$\mbox{lev}_n B$  such that:

(a) the  $x_\alpha^{q^\prime}$,  $\alpha \in v$, are pairwise distinct

(b) for  $\alpha \in u$, $x_\alpha^q \subseteq  x_\alpha^{q^\prime}  $.

This is possible since the  $x_\alpha^q$,  $\alpha \in u$, are distinct in 
$\mbox{lev}_m B$  and  $\vert \mbox{lev}_n B \vert = 2^n = a_n \geq \vert v 
\vert$ 
(cf.\ Definition \ref{de:2.1}(a)). Then define, for  $m \leq i < n$, the
bijection  
$f^{q^\prime} \restriction \mbox{lev}_i B \rightarrow a_i$  such that, if 
$x \in \mbox{lev}_i B$  happens to be  $x_\alpha^{q^\prime} \restriction i$  
for 
some  $\alpha \in u$,  $f^{q^\prime}(x) = t_\alpha^p(i)$. This works since 
the  
$x_\alpha^{q^\prime} \restriction i$, $\alpha \in u$, are distinct 
(recall  $m \leq i$  and the  $x_\alpha^q$,  $\alpha \in u$ 
are distinct), and  the  $t_\alpha^p(i)$,  $\alpha \in u$, 
are distinct by 
Proposition \ref{pr:4.2}.(c).
\endproof

\section{$P^0$,  $P^1$,  and a dense subset of  $P^0 * P^1$}

We define here the forcings  $P^0$  (in  $V$)  and  $P^1$  (in  $V^{P^0}$), 
describe the standard part of  $P^0*P^1$  and find a dense subset  of the 
standard 
part. The central property of the construction is that, on one hand,  
$P^0*P^1$ 
is easily seen to be Cohen and, on the other hand,   $P^0$  and 
$P^1$  adjoin the same generic objects (a function  $\underline f$  and,
for each  $\alpha \in \kappa$, functions  
$\underline t_\alpha  : \omega \rightarrow \omega$  and  
$\underline x_\alpha : \omega \rightarrow 2$) as 
$Q^0$  and  $Q^1$;  this is why  $B(P^0*P^1)$  is isomorphic to  $B(Q^0*Q^1)$.

\begin{definition} \label{de:5.1}   An element of  $P^0$  is a function  
$f$ such that for some  $n \in \omega$ (the height of   $f$,  $\mbox{ht } f$):

(a)  $f$  maps  $\bigcup_{i<n} \mbox{lev}_i B$  into  $\omega$  and for  
$i<n$, 
$f \restriction \mbox{lev}_i B$  is a bijection from  $\mbox{lev}_i B$  onto  
$a_i$.

For  $f$  and  $f^\prime$  in  $P^0$:

(b)  $f \leq f^\prime$ iff   $f \subseteq f^\prime$  (and hence 
$\mbox{ht } f \leq \mbox{ht } f^\prime$).
\end{definition}

Thus for  $K$  $P^0$-generic over  $V$,  ${\underline f}_K = \bigcup K$  is 
a map from the tree  $B$  into  $\omega$, mapping the  $i$\rq th  th level 
of  $B$ 
in a one-one manner onto  $a_i$.  Using the canonical name  $\underline f$ 
for
the generic function  ${\underline f}_K$, we can define the forcing  $P^1$  
in 
$V^{P^0}$.

\begin{definition}  \label{de:5.2} (In  $V^{P^0}$)  $P^1$  is the
finite-support 
product
of the forcings  $P_\alpha^1$,  $\alpha \in \kappa$,  defined as follows. 
An element of  $P_\alpha^1$  is  a pair  
$q_\alpha = ( x_\alpha^{q_\alpha}, t_\alpha^{q_\alpha})$ such that for 
some  $m_\alpha \in \omega$ 
(the height of   $q_\alpha$,  $\mbox{ht } q_\alpha$):

(a)  $x_\alpha^{q_\alpha} \in \mbox{lev}_{n_\alpha} B$ and  
$t_\alpha^{q_\alpha}
\in \mbox{lev}_{n_\alpha}  T$.

For  $q_\alpha$  and  $q_\alpha^\prime$  in  $P_\alpha^1$, 
$q_\alpha \leq q_\alpha^\prime$  iff  

(b)  $\mbox{ht } q_\alpha \leq \mbox{ht } q_\alpha^\prime$, 
$x_\alpha^{q_\alpha} \subseteq x_\alpha^{q_\alpha^\prime}$  and 
$t_\alpha^{q_\alpha} \subseteq t_\alpha^{q_\alpha^\prime}$

(c)  for  $i \in [\mbox{ht } q_\alpha, \mbox{ht } q_\alpha^\prime)$,  
$t_\alpha^{q_\alpha^\prime}(i) = \underline f
(x_\alpha^{q_\alpha^\prime} \restriction i)$.
\end{definition}

\begin{remark} \label{re:5.3} The forcing  $P^0$  is Cohen, since it is
countable 
and
every element has two incompatible extensions. For the same reason, each  
$P_\alpha^1$  is Cohen (in  $V^{P^0}$), hence  $P^1$  is Cohen in  $V^{P^0}$.
It follows that  $P^0*P^1$  is Cohen.              
\end{remark}

\begin{proposition}  \label{pr:5.4} 
The elements of  $\mbox{stp } (P^0*P^1)$  are 
those pairs  
$(f,q)$  such that  $f \in P^0$,  $q = (x_\alpha^q ,t_\alpha^q )_{\alpha \in 
w}$ 
where  $w$  ($= \mbox{dom } q$)  is a finite subset of  $\kappa$, and there 
are natural numbers  $n_\alpha$,  $\alpha \in w$, such that, for  $\alpha 
\in w$:

(a) $\mbox{ht } f \geq n_\alpha$,  $x_\alpha^q \in \mbox{lev}_{n_\alpha} B$, 
and  
$t_\alpha^q \in \mbox{lev}_{n_\alpha} T$ 
(we write  $n_\alpha = \mbox{ht } x_\alpha^q  =  \mbox{ht } t_\alpha^q $).

For  $(f,q)$  and  $(f^\prime, q^\prime)$  in  $\mbox{stp } (P^0*P^1)$,  
$(f,q) \leq (f^\prime, q^\prime)$  iff 

(b) $\mbox{ht } f \leq \mbox{ht } f^\prime$, $\mbox{dom } q 
\subseteq \mbox{dom } q^\prime$, and, for  $\alpha \in \mbox{dom } q$,  
$\mbox{ht } x_\alpha^q \leq \mbox{ht } x_\alpha^{q^\prime}$

(c)  $f \subseteq f^\prime$

(d) for  $\alpha \in \mbox{dom } q$,  $x_\alpha^q \subseteq 
x_\alpha^{q^\prime}$
and  $t_\alpha^q \subseteq t_\alpha^{q^\prime}$

(e) for  $\alpha \in \mbox{dom } q$  and  $i \in [\mbox{ht } x_\alpha^q ,
\mbox{ht } x_\alpha^{q^\prime})$,  $t_\alpha^{q^\prime}(i) = 
f^\prime(x_\alpha^{q^\prime} \restriction i)$.
\end{proposition}

\begin{proposition}  \label{pr:5.5} 
The following subset of  $\mbox{stp } (P^0*P^1)$  
is dense in  $\mbox{stp } (P^0*P^1)$, hence in  $P^0*P^1$.
\[ \begin{array}{rl}
D_P = \{(f,q) \in \mbox{stp} (P^0*P^1) : & \mbox{for all } 
\alpha \in \mbox{dom } q,  
\mbox{ht } x_\alpha^q (= \mbox{ht } t_\alpha^q) = \mbox{ht } f, \\
& \mbox{and the } x_\alpha^q, \alpha \in \mbox{dom } q, 
\mbox{are pairwise distinct}\}. 
\end{array} \]
\end{proposition}

\proof  Let  $(f,q) \in \mbox{stp } (P^0*P^1)$  be given; say with  
$\mbox{ht } f = n$,  $\mbox{dom } q = w$, and  
$\mbox{ht } x_\alpha^q = \mbox{ht } t_\alpha^q = n_\alpha$, for  $\alpha \in 
w$ .
We will find  $(f^\prime, q^\prime) \in D_P$  extending  $(f,q)$
such that  $\mbox{dom } q^\prime = w$  and 
$\mbox{ht } f^\prime = \mbox{ht } x_\alpha^{q^\prime} =
\mbox{ht } t_\alpha^{q^\prime} = N$, where  $N$ is sufficiently large.

To this end, put  $m = \mbox{max } \{  n_\alpha: \alpha \in w \}$  and take 
$N$  so large that  $m \leq N$  and  $\vert w \vert \leq 2^{N-m}$.
Thus we can choose, for  $\alpha \in w$,  $x_\alpha^{q^\prime} \in 
\mbox{lev}_N B$  
such that  $x_\alpha^q \subseteq x_\alpha^{q^\prime}$  and the 
$x_\alpha^{q^\prime}$,  $\alpha \in w$,  are pairwise distinct. 
Fix an arbitrary extension  $f^\prime$  of  $f$  in  $P^0$  of height  $N$.
 Finally define  $t_\alpha^{q^\prime} \supseteq t_\alpha^q$  for  $\alpha 
\in w$ 
 by  $t_\alpha^{q^\prime}(i) = f^\prime (x_\alpha^{q^\prime} \restriction 
i)$, 
 for  $i \in [n_\alpha, N)$.                 \endproof

\section{Conclusion}

According to the sketch of proof given in the introduction, we are left
 with showing that the dense subsets  $D_P$  of  $P^0*P^1$  and  $D_Q$  
 of $Q^0*Q^1$  given in Sections 4 and 5 are isomorphic. This is 
straightforward,
 since  $D_P$  is the following partial order 
(cf.\ Propositions \ref{pr:5.4}, \ref{pr:5.5}). An element  
$\rho$  of
 $D_P$  is, for some finite  $u \subseteq \kappa$  and some  $m \in 
\omega$, 
 a sequence $\rho = (f, (x_\alpha)_{\alpha \in u}, (t_\alpha)_{\alpha \in u}) 
$  
 where

1.  $f$  maps  $\bigcup_{i<m} \mbox{lev}_i B$  into  $\omega$; for  $i<m$, 
$f \restriction \mbox{lev}_i B$  is a bijection from  $\mbox{lev}_i B$  onto  
$a_i$

2.  the  $x_\alpha$,  $\alpha \in u$, are pairwise distinct elements of  
$\mbox{lev}_m B$

3.  the  $t_\alpha$,  $\alpha \in u$, are elements of  
$\mbox{lev}_m T$.

And for  $\rho = (f, (x_\alpha)_{\alpha \in u}, (t_\alpha)_{\alpha \in u}) $  
and
$\sigma = (g, (y_\alpha)_{\alpha \in v}, (s_\alpha)_{\alpha \in v}) $  (with 
domain  $v$  and height  $n$)  in  $D_P$,  $\rho \leq \sigma$  iff
the following hold.

4. $u \subseteq v$  and  $m \leq n$

5. $f \subseteq g$

6. for  $\alpha \in u$,  $x_\alpha \subseteq y_\alpha$  and  
$t_\alpha \subseteq s_\alpha$

7. for  $\alpha \in u$  and  $i \in [m,n)$,  $s_\alpha(i) = 
g(y_\alpha \restriction i)$.

Note that 2. implies that  $\vert u \vert  \leq a_m$; similarly, 
2., 1., and 7.\ imply that the   $s_\alpha$,  $\alpha \in u$,  
are  disjoint above  $m$. Thus, 
up to permutation of coordinates,  $D_Q$  is the same partial order 
(cf.\ Definition \ref{de:4.1} and Propositions \ref{pr:4.2}, \ref{pr:4.3}).


\end{document}